\pdfoutput=1
\documentclass[paper=a4, fontsize=11pt]{article}

\usepackage{comment}
\usepackage[a4paper,pdftex]{geometry}
\usepackage[english]{babel}
\usepackage{listings}
\usepackage{float}
\usepackage{algorithm}
\usepackage{setspace}
\usepackage{indentfirst}
\usepackage{algpseudocode}
\usepackage{amsmath}
 
\usepackage{graphicx}
\usepackage{etoolbox}
\usepackage{multicol}
\usepackage{algorithm,algpseudocode}
\usepackage{amsmath,amssymb}
\usepackage{graphicx}
\usepackage{setspace}
\usepackage{fullpage}

\usepackage{amsthm}
\newtheorem{definition}{Example}
\theoremstyle{definition}
\usepackage{multirow}
    
\usepackage[font=small,labelfont=bf]{caption}
\usepackage{subcaption} 
\usepackage{caption}

\usepackage[
backend=biber,
citestyle=numeric
]{biblatex}
\usepackage{etoolbox}
\patchcmd{\pprintMaketitle}
  {\fi\hrule}
  {\fi\ifvoid\extrainfobox\else\unvbox\extrainfobox\par\vskip10pt\fi\hrule}
  {}{}

\newenvironment{extrainfo}
  {\global\setbox\extrainfobox=\vbox\bgroup\parindent=0pt }
  {\egroup}
\newsavebox\extrainfobox
\title{Importance Sampling based Exploration in Q Learning}

\author{Vijay Kumar$^1$ \and Mort Webster$^2$}
\date{%
    $^1$Industrial and Manufacturing Engineering Department, Pennsylvania State University, Leonhard Building, University Park, Pennsylvania 16802, United States \\%
    $^2$Energy and Mineral Engineering Department, Pennsylvania State University, Hosler Building, University Park, Pennsylvania 16802, United States\\[2ex]%
}

\begin{document}

\maketitle

\begin{abstract}
Approximate Dynamic Programming (ADP) is a methodology to solve multi-stage stochastic optimization problems in multi-dimensional discrete or continuous spaces. ADP approximates the optimal value function by adaptively sampling both action and state space. It provides a tractable approach to very large problems, but can suffer from the exploration-exploitation dilemma. We propose a novel approach for selecting actions using importance sampling weighted by the value function approximation in continuous decision spaces to address this dilemma. An advantage of this approach is it balances exploration and exploitation without any tuning parameters when sampling actions compared to other exploration approaches such as Epsilon Greedy, instead relying only on the approximate value function. We compare the proposed algorithm with other exploration strategies in continuous action space in the context of a multi-stage generation expansion planning problem under uncertainty.
\end{abstract}

\begin{extrainfo}
\textbf{Declarations of interest: none}
\end{extrainfo}

\pagebreak

\section{Introduction}

Multi-stage stochastic optimization refers to the class of problems where the uncertainty evolves dynamically, decisions are adaptive to the realized uncertainty, and there are three or more stages. For many practical situations, multi-stage stochastic optimization may have a multi-dimensional and continuous feasible decision space and either a multi-dimensional continuous uncertainty space or a finite scenario tree with many scenarios.  Such problems will generally not be solvable with exact formulations, such as a deterministic equivalent stochastic program \cite{birge2011introduction}.  A variety of approaches have been developed for solving very large multi-stage stochastic problems in continuous decision space. One class of algorithms builds on the stochastic programming paradigm and extends two-stage decomposition schemes such as Benders Decomposition \cite{van1969shaped} or Progressive Hedging \cite{rockafellar1991scenarios}. However, these methods depend on the existence of dual variables in the subsequent stage or on efficient and stable quadratic programming solvers and can be difficult to extend to problems with non-convex constraints, although some recent work has been done to address this \cite{zou2018partially,zou2019stochastic,liu2017multistage}. 

Another class of methods applied to continuous state and decision spaces has developed from Stochastic Dynamic Programming and Markov Decision Problems \cite{bellman2015applied}. Two approaches for solving dynamic programming problems in continuous spaces are discretizing the state and action space~\cite{puterman2014markov} to obtain the policy, and use of analytical closed form solutions for the optimal policy, when such a solution exists~\cite{lewis2012optimal}. The closed form solution can be applied to specific problems such as Linear Quadratic Control problems whose analytical solution can be found by solving the Riccati equations. However extension of this approach to problems with a non linear structure is difficult because an algebraic solution cannot be obtained. An advantage of discretization is that it does not depend on the problem characteristics (e.g. convexity of the space), but it suffers from the Curse of Dimensionality and  can become computationally intractable for problems with high dimensional continuous decision and state spaces. For very large problems with both continuous state and decision spaces, Monte Carlo sampling-based approaches such as Approximate Dynamic Programming (ADP) methods \cite{powell2007approximate} or Reinforcement Learning (RL) methods \cite{sutton2018reinforcement} can find good approximate solutions to the original problem by sampling uncertainties and decisions and using low-order functional approximations to the value or recourse function. These methods are generally computationally less expensive than an exact Stochastic Dynamic Programming approach and can be used in a multi-stage adaptive framework.  

The computational advantage of ADP lies partly in using sampled costs for specific states and actions to estimate a functional approximation that extrapolates the information to other potential states and actions not sampled.  Functions can be used to approximate either the value function, which maps states to optimal expected values, or a state-action value function, which maps a state-action pair to the expected value of choosing that action in that state and subsequently following the optimal policy.  We focus on state-action approximations applied to continuous state and action spaces within a Q-Learning framework \cite{watkins1992q}.  

Q-learning uses a low dimensional approximation of the objective function and progressively updates the approximation as it observes the costs for sampled states and actions The goal of the algorithm is to learn a policy that mimics the optimal policy obtained by exact Dynamic Programming but with less computational effort. When the complete enumeration of the feasible decision space is computationally prohibitive, an efficient sampling methodology is required. Q-learning algorithms must balance the trade-off between sampling state-action pairs where the variance is large to improve the quality of the estimate and exploiting the existing estimate to minimize the cost. This is commonly known as the dilemma of exploration vs exploitation~\cite{powell2007approximate,kaelbling1996reinforcement,bertsekas2012approximate}. Several approaches have been developed to navigate this tradeoff, of which one widely used method is Epsilon Greedy \cite{watkins1989learning}.  Epsilon greedy randomizes the strategy between pure exploration, i.e., sampling any feasible action, and pure exploitation, i.e., choosing the optimal action from the current estimated state-action value function. The relative ratio of explore vs. exploit samples will generally be a function of the iteration count.  One challenge for epsilon greedy and similar methods is that in early iterations, the value function approximation is poor, so that the exploit samples may not be useful or can be misleading, and in later iterations as the approximation improves, exploration samples do not utilize the global or past information and may inefficiently continue to sample sub-optimal actions.  A second challenge is that these methods each have algorithmic parameters that must be heuristically tuned for each application, and convergence can be quite sensitive in practice to the values chosen. A third challenge to epsilon greedy is that an additional optimization problem must be solved to obtain each exploit sample action, which can increase the computational cost considerably if the value function approximation is non-linear or non-convex.

Importance Sampling (IS) \cite{rosenbluth1955monte,hammersley1954poor} is a variance reduction technique that disproportionately samples from a distribution other than the "true" distribution to estimate the expected value with fewer samples. Importance Sampling has been widely applied in Stochastic Programming \cite{infanger1992monte,parpas2015importance} and Reinforcement Learning \cite{precup2000eligibility,precup2001off} to provide variance reduction. IS within Stochastic Programming has been applied to efficiently sample large uncertainty spaces, while in these methods the decisions are obtained by solving linear programs. In Reinforcement Learning, Importance Sampling has been used to eliminate bias when the behavioral or sampling policy is different from the optimal policy. In the update step to improve the approximate value function, the sampled rewards are weighted using the density ratio to reflect the likelihood of sampling the state-action pair. However, these applications use other methods for selecting the samples.  Importance Sampling offers a potentially powerful means of sampling actions, making use of information learned in prior iterations without going to the extreme of pure exploitation. The approximate value function can be used within an accept-reject framework to explore actions disproportionately in the neighborhood of the current optimum, with some exploration of other actions still occurring.

In this paper, we present a novel sampling scheme to select actions in continuous space within a Q-learning algorithm based on importance sampling.  We propose to use an Importance Sampling weighted accept-reject method to sample actions; i.e., importance sampling defines the sampling policy rather than weighting observed rewards from some other sampling policy. We exploit the characteristic of disproportionate sampling in IS to balance exploration and exploitation. We apply this method to multi-stage Stochastic Generation Expansion Planning (GEP), the problem of choosing electricity generation  investments in each stage under uncertainty. GEP provides an example of a problem with multi-dimensional continuous feasible action and state spaces.  In the context of this application, we compare the performance of our proposed sampling scheme to two alternative methods used in Approximate Dynamic Programming and Reinforcement Learning for continuous decision spaces, epsilon greedy \cite{watkins1989learning} and epsilon decay.  In particular, we demonstrate that the latter methods for decision selection can sometimes result in convergence to sub-optimal policies, and that different algorithmic parameter choices can result in widely differing policies from replications of the same problem. In contrast, we show that the importance sampling-based approach converges to higher quality approximate solutions with less variance across repetitions, and that the algorithm increases the density of sampling near-optimal actions as the value function approximation improves.

The contributions of this paper are:
\begin{itemize}
    \item A novel algorithm based on importance sampling to sample actions in continuous decision space for Q-learning.
    \item Application of the method to multi-stage stochastic generation expansion planning and comparison with epsilon greedy sampling methods.
    \item Demonstration of convergence to more accurate approximate solutions with less variance, without any algorithmic parameters.
\end{itemize}

The remainder of the paper proceeds as follows. Section 2 reviews sampling methods within Approximate Dynamic Programming and Reinforcement Learning algorithms and the use of Importance Sampling in Reinforcement Learning. Section 3 introduces Markov Decision Processes.  Importance Sampling is defined in Section 4. Section 5 presents our Q-learning Importance Sampling method. The application problem, multi-stage generation expansion planning, is introduced in Section 6. The proposed algorithm is compared with other sampling methods in Section 7 and Section 8 provides a concluding discussion.

\section{Literature Review}

Approximation architectures provide a computationally tractable approach to solve problems with high dimensional state-action spaces when exact algorithms cannot provide solutions in reasonable computational time. Approximate Dynamic Programming \cite{powell2007approximate} and Reinforcement Learning \cite{sutton2018reinforcement} offer a general framework that can provide solutions to complex, high dimensional optimization problems under uncertainty. The algorithms within this class of methods use a combination of statistical methods for sampling and estimation, optimization methods, and lower dimensional representation of the objective function to obtain approximate solutions to the original problem.  

The performance of Monte Carlo-based methods depend critically on efficiently adaptively sampling the state-action space to learn the best approximation of the true value function. The challenge for any sampling policy is that actions that appear sub-optimal at the beginning of the algorithm when the approximation is less accurate may in fact be optimal and vice versa. It is important to sample actions not only based on the current expected reward, but to balance this with exploring the action space to improve the policy. This trade-off is known as the ‘explore vs exploit’ problem and has received considerable attention in the research community \cite{powell2007approximate,kaelbling1996reinforcement,bertsekas2012approximate}. Many sampling methods and heuristics have been developed in the literature to better balance exploration and exploitation, including epsilon-greedy, R-Max \cite{brafman2002r} and E$^{3}$\cite{kearns2002near}.  R-Max and E$^{3}$ require enumeration of the state space, which is possible for some reinforcement learning applications, but becomes computationally intractable for large or continuous state spaces. Epsilon greedy uses an algorithmic parameter to prescribe the relative frequency of exploration and exploitation.  Epsilon decay prescribes a trajectory for the explore/exploit ratio that changes with the iteration count. In practice, the algorithmic parameters in these methods that specify the relative frequency of exploration and the rate of change in the ratio can be difficult to tune for each application. Epsilon greedy methods have an additional disadvantage for some applications; every exploit sample is generated by solving for the optimal action relative to the current value function approximation. If the approximation architecture is non-linear, this can increase the computational burden. 

An active area of research in the Reinforcement Learning literature to promote exploration is the application of an exploration bonus \cite{gehring2013smart, strehl2008analysis} or intrinsic motivation \cite{singh2010intrinsically,bellemare2016unifying} to assist the algorithm in exploring areas which would not otherwise be visited but that will improve the approximation and future decisions. The fundamental idea in these algorithms is to add a bonus term to the Bellman Optimality equation based on the state-action visit counts. Although theoretical guarantees exist for tabular implementations \cite{strehl2008analysis}, only recently has it been applied to high dimensional state spaces by using a density model over the state space \cite{martin2017count,ostrovski2017count,bellemare2016unifying}, a hashing method to cluster similar states \cite{tang2017exploration} or by definition of another functional approximation for the exploration count \cite{gehring2013smart}. All these methods require the state-action space to be countable, limiting the applicability to continuous settings. Further, approaches that rely on modifying the reward as an incentive to explore may not work well for problems in which actions strongly impact future rewards and decisions, such as the capacity expansion problem. 

Another strategy for exploration has been the use of Bayesian methods in Approximate Dynamic Programming.  Early examples defined a prior over the transition probabilities \cite{silver1963markovian} and applied this approach to multi-armed bandit problems \cite{duff1997local}. Variants of this methodology also define a prior over the value function itself \cite{dearden1998bayesian}. Bayesian methods were extended to account for correlation in the impact of individual actions across many states \cite{ryzhov2010approximate}. The Knowledge Gradient metric was developed to sample actions within a Bayesian framework \cite{ryzhov2019bayesian}. Although the evaluation of the Knowledge Gradient metric does not require new parameters, the Bayesian model to which it is applied requires the definition of a prior over the value function and also requires the specification of a covariance matrix, which requires practical knowledge of the underlying system. The Knowledge Gradient approach can be extended to continuous state spaces, but requires the enumeration of the action space and cannot be applied to continuous action spaces without discretization.

Reinforcement Learning for problems with continuous action spaces has applied policy gradient methods~\cite{sutton2000policy,silver2014deterministic} and cross entropy search methods~\cite{simmons2019q,kalashnikov2018qt}. However these methods have faced challenges when defining a policy over multi-dimensional action spaces or within a constrained action space. Actor-critic~\cite{lillicrap2015continuous,haarnoja2018soft} methods integrate value function approximation and policy gradient algorithms for such multidimensional action spaces, but the iterative estimation of both value and policy function approximations can substantially increase the computational demands. Recently, new variants of value function approximation methods that use neural network architectures have sampled actions using normalized advantage functions~\cite{gu2016continuous}, convex inputs to a neural network~\cite{amos2017input} or mixed integer-based formulations~\cite{ryu2019caql}. All the above require the optimization of an approximate analytical function to generate proposed samples during the exploit phase and can result in sub-optimal actions if the feasible action space does not satisfy the assumed structure of the analytical function. The above methods are specific to neural networks, and may be difficult to adapt to other approximation architectures.

Importance Sampling has been widely applied in the Reinforcement Learning community since it was first introduced using per-decision importance sampling in tabular reinforcement learning \cite{precup2000eligibility} and subsequently expanded to include linear function approximations \cite{precup2001off}. The problems addressed in the reinforcement learning community typically have long episodes or sample paths; in this context, per-decision importance sampling provides an unbiased estimator with variance reduction so that fewer samples are required to converge \cite{precup2000eligibility}. Early success with this approach motivated the development of variance reduction techniques such as weighted importance sampling \cite{mahmood2014weighted}, adaptive importance sampling \cite{hachiya2009adaptive}, emphatic TD learning \cite{sutton2016emphatic}, doubly robust methods \cite{jiang2016doubly}, conditional importance sampling \cite{rowland2020conditional} and incremental importance sampling \cite{guo2017using}. All these methods are focused on facilitating an unbiased update step when the sampled reward is the result of an action chosen from a policy other than the current optimal policy.  Prediction or behavior policies, which generate the sample state-action pairs, have included epsilon greedy \cite{guo2017using}, uniform random sampling \cite{jiang2016doubly,mahmood2014weighted,xie2018marginalized,gelada2019off}, estimation through multiple trajectories \cite{hanna2019importance}, softmax distributions \cite{xie2019towards} or specified probability distributions of actions \cite{precup2000eligibility,precup2001off,sutton2016emphatic,mahmood2017multi}.  These methods tend to sample actions from the full feasible space without conditioning on the previous samples observed, except when using the current approximate value function to exploit. For example, if we consider epsilon greedy as the behavior policy, an explore sample may randomly select any feasible action. But this sampling policy provides no additional guidance on the relative benefit between actions for improving the approximation. In large or continuous action spaces when every action cannot be visited, the behavior policy should sample the actions that will provide the maximum improvement to the value function approximation. In this paper, we demonstrate the use of importance sampling to guide the sampling of actions in a continuous space, rather than simply correcting the bias for sample rewards obtained by some other means.

\section{Markov Decision Processes}
Markov Decision problems are commonly used for modelling multi-stage optimization problems under uncertainty. Markov Decision Processes are defined by the following elements:

\begin{itemize}
    \item Stages ($t \in T$): Specify the specific points in time where an action is chosen by the agent based on the information available up to that point.
    \item State ($s \in S$): The vector of information that makes the process Markovian; i.e., it contains all information required to specify the reward to be received and the state of the process in the next stage that will result from each feasible action.  
    \item Action ($a \in A$):  A feasible action that can be chosen by the agent in the current state and stage.
    \item Exogenous Information ($\omega \in \Omega$): The new information that is obtained between stages $t-1$ and $t$ for every state; represents the uncertainty in the system. 
    \item State Transition Function: Represents how the state variable changes as a result of the chosen action and the realized exogenous information. Often takes the form of a probability distribution $P_{s,s'}(s,a,\omega)$.
    \item Cost or Reward ($r$): This is the immediate reward or cost to the agent as a result of choosing action $a$ in state $s$; may be deterministic or stochastic. 
    \item Policy ($\pi$): A policy maps states to a probability distribution over its actions. The mapping can be invariant over time (stationary) or can change as a function of its history (non-stationary).
\end{itemize}

There are several distinct classes of MDPs, depending on whether it is a finite horizon or an infinite horizon problem and on the form of the total cost over all stages \cite{puterman2014markov}.   Here, we restrict our focus to finite horizon problems in which the total cost is the discounted sum of costs in each stage.  We will therefore include the subscript $t$ to indicate the dependence on stage with the probability transition function defined as $P_{s_{t},s_{t+1}}$ for the state transition from stage $t$ to $t+1$. We implicitly include the uncertainty $\omega_{t}$ within the state variable $s_{t}$ and the transition function $P_{s_{t},s_{t+1}}$ to simplify notation. The formulation applies to either total reward maximization or total cost minimization; without loss of generality we present the cost minimization notation to be consistent with the application in this paper. The expected value of following any feasible policy $\pi$ is

\begin{equation}\label{Value function}
v^{\pi}_{t}(s) = E_{\pi}\Bigg\{\sum_{t=t^\prime}^{T} \gamma^{t}r_{t}\Bigg | s_{t}=s \Bigg\}
\end{equation}

This represents the expected value of the sum of discounted costs when the system begins in state $s$ and then selects actions based on $\pi$ for all remaining stages. In Equation (\ref{Value function}), $r_{t}$  is the cost and $s_{t}$ is the state at time period $t$ and  $\gamma^{t}$ is the discount factor. 

The objective is to solve for the optimal policy $\pi^{*}$ that minimizes the value function (\ref{Value function}):
\begin{equation}\label{Optimal Value function}
v^{*}_{t}(s) = \min_{\pi}v^{\pi}_{t}(s) \ \ \ \ \textrm{for} \ \ \textrm{all}\ \ \ s\in S, t \in T
\end{equation}

The exact solution of the finite horizon MDP in (\ref{Optimal Value function}) can be solved with dynamic programming, a stage-wise decomposition of the full problem into a sequence of smaller problems for each stage and state.  The solution for each subproblem is found via the Bellman Equation with $P_{s_{t},s_{t+1}}$ representing the probability distribution of the state transition function:

\begin{equation}\label{Bellman Equation}
v^{*}_{t}(s_t) = \min_{a_t}E_{P_{s_{t},s_{t+1}}}\Bigg\{ r_{t}(s_t,a_t) +  \gamma v^{*}_{t+1}(s_{t+1}) \Bigg\} \ \ \ \ \textrm{for} \ \ \textrm{all}\ \ \ s\in S, t \in T
\end{equation}

The exact solution of (\ref{Bellman Equation}) can be obtained by backward induction, starting in the final stage $T$ and solving for every state, and then moving back one stage and repeating until solved for $t=1$.   This approach is not computationally tractable once the number of stages, states, actions, and exogenous information signals become too large. Approximate algorithms (e.g., ADP, RL) focus on iteratively sampling state-action pairs $(s_t,a_t)$, observing the sampled costs $r_t$, and updating a low-order approximation of the value function $\hat{v}_t$.   

An alternative form of the value function known as the state-action value function, is useful for some applications:

\begin{equation}\label{Bellman Q}
q_{t}(s_t,a_t) = E_{P_{s_{t},s_{t+1}}} \{ r_{t}(s_t,a_t) + \gamma \min_{a_{t+1}}q_{t+1}(s_{t+1},a_{t+1}) \} \ \ \ \ \textrm{for} \ \ \textrm{all}\ \ \  t \in T
\end{equation}

The solution minimizes

\begin{equation}\label{Optimal State Action Value function}
v_t^{*} = q^{*}_{t}(s,a) = \min_{\pi}q^{\pi}_{t}(s,a) \ \ \ \ \textrm{for} \ \ \textrm{all}\ \ \ s \in S, a \in A, t \in T
\end{equation}

A sampling-based approach approximates the state-action value function $q(s,a)$ with a low-order projection in terms of a set of features/basis functions that depend on the state and action~\cite{bertsekas2012approximate}. The approximations can be linear or non-linear and an extensive list of functional forms is described in Sutton and Barto~\cite{sutton2018reinforcement}. The simplest functional form is a linear approximation with coefficients or weights $\theta$.  The algorithm solves for the approximate function $\hat{q}$: 
\begin{equation}\label{approximation}
 \hat{q}_{t}(s,a,\theta) \approx q_{t}(s,a) \ \ \ \ \textrm{for} \ \ \textrm{all}\ \ \  t \in T
\end{equation}

The approximation is updated iteratively by sampling specific state-action pairs with temporal difference learning~\cite{sutton2018reinforcement} as in Equation (\ref{Update Q Value}),
\begin{equation}\label{Update Q Value}
 \hat{q}_t(s,a,\theta) \leftarrow (1-\lambda)\hat{q}_t(s,a,\theta) + \lambda(r_t(s_t,a_t) + \gamma\min_{a_{t+1}}\hat{q}_{t+1}(s_{t+1},a_{t+1},\theta)) \ \ \ \ \textrm{for} \ \ \textrm{all}\ \ \  t \in T  
\end{equation}

The second term in Equation \ref{Update Q Value} can be estimated by looking $n$ steps ahead, where $n$ can take any value from 1 (Equation \ref{Update Q Value}) to $T$ (Equation \ref{Update Q Value 2}). 

\begin{equation}\label{Update Q Value 2}
 \hat{q}_t(s,a,\theta) \leftarrow (1-\lambda)\hat{q}_t(s,a,\theta) + \lambda(r_t(s_t,a_t) + \gamma r_{t+1}(s_{t+1},a_{t+1}) ... r_{T}(s_{T},a_{T}))  \ \ \ \ \textrm{for} \ \ \textrm{all}\ \ \  t \in T 
\end{equation}

In equation (\ref{Update Q Value}, \ref{Update Q Value 2}), we can either use a single state-action sample for the update or we can perform an expected update using a batch of samples. Although the computational complexity is higher for the batch sample update, performance improvement is often observed. 

\section{Importance Sampling}
Consider the general problem of estimating the expected value of a function $h(x)$ of a random variable $x$ with probability density function $f(x)$.

\begin{equation}\label{IS1}
\mathcal{L}=E_{f}[h(X)]=\int h(x)f(x)dx	
\end{equation}

For estimating extreme values or rare events, the number of samples of $x$ required may be very large.   In these cases, a useful strategy is to sample from a different distribution $g(x)$ such that $h(x)f(x)$ is dominated by $g(x)$.  Because the expectation is taken with respect to the original distribution $f(x)$, each sample must be weighted by its relative likelihood:

\begin{equation}\label{IS2}
\mathcal{L} = \int h(x) \dfrac{f(x)}{g(x)}g(x)dx	
\end{equation}

Then an unbiased Monte Carlo estimator of $\mathcal{L}$ is

\begin{equation}\label{IS3}
\hat{\mathcal{L}} = \sum_{i} { h(x_{i}) \dfrac{f(x_{i})}{g(x_{i})}g(x_{i}) }	
\end{equation}

The choice of the importance density function $g(x)$ is critical in determining the amount of variance reduction that can be achieved.   It can be shown \cite{kroese2011handbookmc} that the variance of $\hat{\mathcal{L}}$ is minimized by sampling from

\begin{equation}\label{IS4}
g^*(x) = \dfrac{\lvert h(x) \rvert f(x)}{\int \lvert h(x) \rvert f(x) dx }	
\end{equation}

The minimum variance density has long been elusive because of the curse of circularity: the denominator in (\ref{IS4}) is the very quantity we are estimating.  Various heuristics have been used to select "good" importance sampling densities (see \cite{kroese2011handbookmc}).  More recently, it has been shown that Markov Chain Monte Carlo methods can be used to sample from the optimal density function \cite{parpas2015importance} up to a proportionality constant.  Note that if the original density $f(x)$ in (\ref{IS4}) is uniform (i.e., $f(x)$ is a constant), than the optimal density function $g(x)$ is simply a normalized version of the function $h(x)$ to be estimated.  We exploit this fact below in our algorithm.

In addition to statistical estimation and simulation, importance sampling has been used extensively in reinforcement learning \cite{sutton2018reinforcement}.  In many RL algorithms, the expected value (the second term in the Bellman equation) is estimated by looking $n$-steps ahead, as shown in Equation \ref{Update Q Value 2}. When temporal difference learning is applied for two or more steps ahead, additional state-action pairs must be sampled using the current approximation of the optimal policy. Efficiency is improved by sampling state-action pairs for the additional steps from a policy different from the optimal policy, known as an off-policy approach.   Because these samples are used to estimate the expected value of future decisions, the likelihood of future rewards must match the likelihood of visiting those states in the learned/optimal policy. Therefore, in the prediction step of these algorithms, off-policy sampling uses the importance sampling ratio to correct these estimates~\cite{mahmood2014weighted,precup2001off}.

The original method of importance sampling in estimation (\ref{IS2}-\ref{IS3}) performs two distinct tasks. The first is to \emph{sample} more efficiently by focusing the measurement or computation on the most useful regions of the sample space.   The second is to then \emph{correct the bias} from disproportionately sampling in the estimated expected value.  The applications of importance sampling in off-policy reinforcement learning only perform the latter of these functions, while using other methods (e.g., epsilon greedy or uniformly random) for selecting the state-action samples~\cite{guo2017using,mahmood2014weighted,gelada2019off}.

In the proposed method, we focus on the former objective, efficiently sampling actions.  Building on the insight that the minimum variance sampling density $g(a)$ is at least proportional to the function being estimated $\hat{q}(s,a)$, and adapting the MCMC sampling approach from \cite{parpas2015importance}, we use an accept-reject algorithm to adaptively sample actions using the current estimated approximate Q-function $\hat{q}$.

We apply this sampling approach to a one-step Q-learning framework.  As noted in \cite{sutton2018reinforcement} (see Sections 7.3 and 12.10), when sampling off-policy in one-step algorithms, sample action $a$ only impacts the reward $r(s,a)$ and the state transition $P_{s,s'}(s,a)$ and therefore the importance sampling ratio is not needed to correct the bias. In other words, in a one-step Q-learning framework, we are estimating the value of $\hat{q}(s,a)$ for the action $a$ that has been sampled already and not the expected value of many different values of $a$. Further, Q-learning is directly approximating the optimal Q-function($q$) and is sampling over the state/uncertainty transitions, once the action $a$ is selected using the off policy approach. Therefore, we do not need to perform the bias correction step with the importance sampling ratio.  Note that the correction would be required if i) one also applied importance sampling to the exogenous uncertainty that determines the state transitions, or ii) if applying this sampling approach to other RL algorithms with two or more steps ahead sampled off-policy.  In these cases, one could adapt the Kernel Density Estimation procedure presented in \cite{parpas2015importance}.

\section{Q-learning Importance Sampling}
We present the Q-learning Importance Sampling (QIS) method (Algorithm 1). The name of the algorithm derives from the fact that it uses the function approximation $\hat{q}$ within an importance sampling framework to sample actions. The QIS algorithm uses accept-reject sampling, where the likelihood of accepting a random action given the current state is proportional to the current approximation of the Q-function $\hat{q}$. Proposed sample actions $\xi$ are generated, and the approximate value  $\hat{q}(\xi)$ for the proposed action is obtained using the current approximation $\hat{q}$. Proposed samples are accepted or rejected based on the $qRatio$, as defined in Equation \ref{IS eqn 2},

\begin{equation}\label{IS eqn 2}
     qRatio =\dfrac{qMax-\hat{q}(\xi)}{qMax-qMin} 
\end{equation}

\noindent{where $qMax$ and $qMin$ are the minimum and maximum approximate $\hat{q}$ function values determined each iteration from all samples observed so far.
In early iterations when the function approximation is poor, the estimates of $qMin$ will also be poor and the probability of a sample in the entire space being accepted will be high because the denominator of $qRatio$ would be small. As more samples are obtained, the estimate of $qMin$  will improve and the probability of proposed samples with values closer to $qMin$ will become more likely to be accepted than samples with values closer to $qMax$, although no proposals will have zero probability of being accepted. The result is a shift in the distribution of sampled actions from more uniformly distributed in the early iterations when the approximation is poor to distributions with higher density near the optimal action and lower density elsewhere in later iterations.  Unlike epsilon greedy, this method will never perform a pure exploit sample, sampling exclusively the optimal action for a state as determined from the approximate Q-function $\hat{q}$.}

\begin{figure}[H]
\centering
\begin{subfigure}{.475\textwidth}
  \centering
  \includegraphics[width=\textwidth]{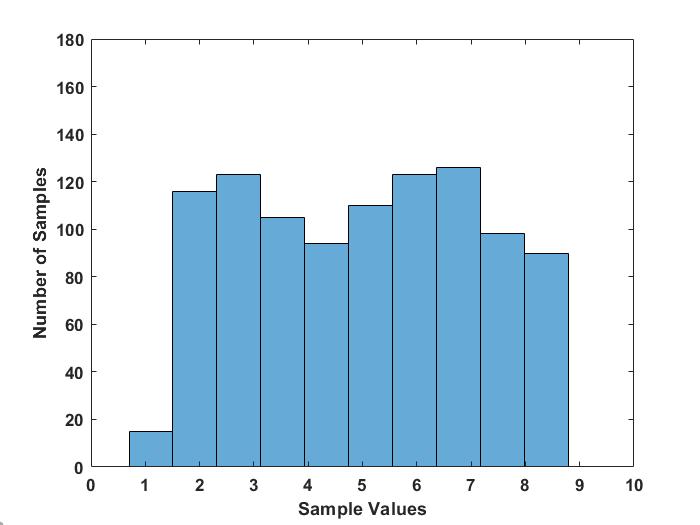}\hfill
  \caption{Initial 1000 Samples}
  \label{fig: Example_First}
\end{subfigure}%
\begin{subfigure}{.475\textwidth}
  \centering
  \includegraphics[width=\textwidth]{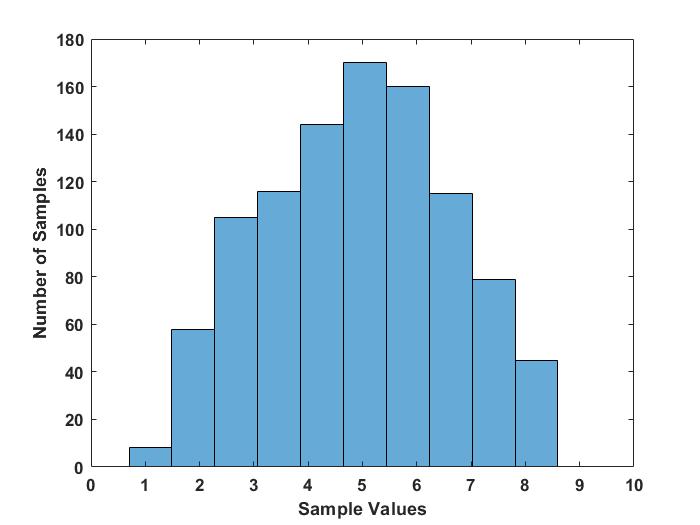}\hfill
  \caption{Final 1000 Samples}
  \label{fig: Example_Last }
\end{subfigure}
\caption{Initial and Final 1000 Samples of QIS algorithm applied to Example \ref{Example-QIS}.}
\label{fig: Example Samples}
\end{figure}

We illustrate this characteristic of the QIS algorithm on a small one dimensional problem in Example \ref{Example-QIS}. Figure \ref{fig: Example Samples} shows the initial and final 1000 samples accepted in Example \ref{Example-QIS} for one representative replication. We observe that the initial 1000 accepted samples are approximately uniformly distributed and the distribution of the final 1000 accepted samples exhibits a higher concentration of samples around the optimal action.

\begin{definition}\label{Example-QIS}
        Minimize the following function: 
        \begin{equation}\label{Example Equation 1}
            Q(x) = 25+(x-5)^2 
        \end{equation}\
Assume that the initial estimates of the minimum and maximum Q-function values are equal to 35 and 40, respectively. The true minimum of (\ref{Example Equation 1}) is 25 at x=5. We perform K=5 iterations of QIS with M=1000 samples at every iteration. Figure \ref{fig: Example Samples} shows the initial and final 1000 samples accepted by QIS for one representative replication. 
\end{definition}

The detailed steps of QIS are given in Algorithm 1.  The algorithm iteratively samples each stage $t$ within each iteration $k$. Steps 1-2 of the algorithm initialize the parameter values for the algorithm, the Q-function approximation $\hat{q_{t}}$ coefficients for stages $t = 1..T$, the approximate minimum ($qMin_{t}$) and maximum ($qMax_{t}$) function values for each stage, and the initial state of the system $s_{0}$. If bootstrapping is not performed, the approximations for all stages will have coefficients of zero and approximate minimum($qMin_{t}$) and maximum($qMax_{t}$) function values assigned the values 0 and 1, respectively.

\begin{algorithm} 
  \begin{algorithmic}[1]
    \State Initialize $\hat{q_{t}}$ for $t=1..T$
    \State Initialize $qMin_t$ and $qMax_t$ for $t=1..T$ and $s_{0}$
    \For{$k=1..K$}
     \State \textbf{Begin Forward Pass}
        \For {$t=1..T$}
        \If{$t = 1$}
         \For{$m = 1..M$} 
        \State $s_{t}^{m} = s_{0}$
        \EndFor
        \Else
        \For{$m = 1..M$}
        \State Sample $\omega_{t}^{m} \in \Omega_{t}$
        \State $s_{t}^{m} = [a_{t-1}^{m} \ \ \omega_{t}^{m}$ ]
        \EndFor
        \EndIf
            \While{$m \leq M$}
                \State Propose an action $\xi$ and calculate approximate value $\hat{q}(\xi)$ using $\hat{q}_{t}$.
                \If{$\hat{q}(\xi) > qMax_{t}$}
                    \State $qMax_{t} = \hat{q}(\xi)$
                \EndIf
                \If{$\hat{q}(\xi) < qMin_{t}$}
                   \State $qMin_{t} = \hat{q}(\xi)$
                \EndIf
                \State Calculate $qRatio =\dfrac{qMax_{t} - \hat{q}(\xi)}{qMax_{t} - qMin_{t}}$
                \State Generate a random number $U$ $\sim$ U(0,1)
                     \If{$qRatio > U$}
                        \State $m = m + 1$
                        \State $a_{t}^{m} = \xi$
                        \State Calculate $r_{t}^{m}(s_{t}^{m},a_{t}^{m})$. 
                        \Else
                        \State Proposed action $\xi$ is rejected.
                    \EndIf
            \EndWhile    
        \EndFor
    
    \State \textbf{End Forward Pass}
    \State \textbf{Begin Backward Pass}
    \For{$t=T..1$}
    \State Update $\hat{q_{t}}$ using M samples in forward pass using equation \ref{Update Q Value}.
    \State Obtain new estimates of $qMax_{t}$ and $qMin_{t}$ using all samples for iterations 1 to $k$.
    \EndFor
    \State \textbf{End Backward Pass}
    \EndFor
  \end{algorithmic} 
  \caption{\textbf{QIS}}
  \label{alg:algorithm1}
\end{algorithm}

For each iteration $k$ and stage $t$, the algorithm performs a forward pass (Steps 4-35) in which it selects $M$ samples of state-action pairs for each stage, followed by a backward pass (Steps 36-41) in which it uses the costs obtained for these samples to update $\hat{q_{t}}$. In the forward pass, for all stages $t \geq 2$, $M$ samples of the exogenous uncertainties $\omega_{t} \in \Omega_{t}$ are generated to create the state vector $s_{t}^{m}$ (Steps 6-15). Steps 16-33 generate sample actions for each sample state. An action $\xi$ is proposed and its value, $\hat{q}(\xi)$, is calculated using the current value function approximation $\hat{q}_{t}$. $qMax_{t}$ and $qMin_{t}$ are updated if $\hat{q}(\xi)$ is greater than $qMax_{t}$ or less than $qMin_{t}$. The proposed sample action is accepted if its $qRatio$ (\ref{IS eqn 2}) is greater than a sample from $U(0,1)$ (Steps 24-30). The cost for this state-action pair is $r_{t}^{m}(s_{t}^{m},a_{t}^{m})$. The process iterates until $M$ samples are accepted for each stage. Steps 36-41 constitute the backward pass of the algorithm, in which the function approximation $\hat{q_{t}}$ is updated using temporal difference learning (\ref{Update Q Value}) with the $M$ sample costs obtained in the forward pass. New estimates for $qMax_{t}$ and $qMin_{t}$ are obtained using the updated approximation $\hat{q_{t}}$. The algorithm repeats until the maximum number of iterations $K$ is reached. 

The QIS algorithm has two components: 1) the accept-reject scheme using $qRatio$ in the forward pass to obtain new sample actions, and 2) in the backward pass, after the value function approximation, $\hat{q_{t}}$, is updated, the estimates of $qMax_{t}$ and $qMin_{t}$ are updated using all the samples observed from previous iterations. This update of $qMax_{t}$ and $qMin_{t}$ necessarily entails re-evaluating the new approximate Q function values for every prior state-action sample, because these values will change for the updated approximation, which will alter the approximate maximum and minimum observed Q-function values for each stage. However, this update procedure becomes computationally expensive as the number of iterations and cumulative samples increase. To reduce the computational burden, we propose a variant of QIS that updates the estimates of $qMax_{t}$ and $qMin_{t}$ every $\hat{K}$ iterations, where $\hat{K}$ can take values between $1$ (QIS) and $K$ (never update). The required computation for re-evaluation is decreased by a factor of at least $\hat{K}$. We call this algorithm QIS-RE because of the reduced frequency of updating $qMax_{t}$ and $qMin_{t}$.

There are three advantages to the QIS algorithm.  The first advantage is that it can be applied to continuous action spaces because it does not require complete enumeration or discretization of the feasible actions. A second advantage is that it does not require the solution of an optimization problem during the sampling step. In the epsilon greedy algorithm, when an action is to be sampled using exploitation, the optimal action for the current state must be determined using the current value function approximation $\hat{q}$.  For many applications, the functional form of the Q-function approximation, $\hat{q}$, will include non-linear terms to accurately represent the true Q-function, $q$; in such cases, a non-linear program must be solved to obtain the next sample action. In QIS, the current value function approximation, $\hat{q}$, is evaluated to determine the $qRatio$ for a proposed sample, but it is not necessary to solve an optimization problem.  A third advantage to the QIS algorithm is the absence of  tuning parameters for sampling actions. In epsilon greedy, $\epsilon$ is a heuristic parameter which is difficult to tune for each specific problem. ADP/RL algorithms require other assumptions, for example the basis function to be used for the low dimensional approximation of the true Q-function, $q$ or the parameter $\lambda$ in temporal difference updating (\ref{Update Q Value}). Each additional heuristic algorithmic parameter further increases the difficulty in achieving convergence and good approximate solutions for each new application.

\section{Application: Multi-stage Generation Expansion}
We demonstrate the QIS algorithm with an application to a multi-stage stochastic electricity generation expansion problem.  Stochastic Generation Expansion Planning (SGEP) is a long-term investment problem that solves for the optimal new capacity investments in each stage that minimize the expected total (capital plus operating) cost with respect to uncertainty in some parameters. The decision vector in each stage is the amount of new capacity to build of each generation type. We first present a mathematical formulation of SGEP, and then describe a reformulation of this problem as a Markov Decision Process.

The objective function (\ref{Objective}) minimizes the sum of expected investment and operations costs over all investment periods $t \in T$, all generation technology types $g \in G$, and all scenarios $\sigma \in \Sigma$. The operations cost is summed over demand load blocks $l \in L$, with $H_{l}$ hours in demand load block $l$. Scenarios are used to represent the uncertainty, and each scenario is associated with a probability $p(\sigma)$. The decision $y_{t\sigma}^g$ represents the new capacity added for each generation type $g$ in scenario $\sigma$ in investment period $t$. Hourly demand within each year is represented by a load duration curve, in which electricity demand $D_{lt}$ is specified for each demand block $l \in L$, and the demand blocks represent hours with different levels of demand within the year. We assume a 2 \% per annum growth rate in aggregate annual demand. The long-term investment problem includes the dispatch sub-problem, choosing $x_{lt\sigma}^g$, the energy output from generation technology type $g$ in scenario $\sigma$ to meet the demand in block $l$ in investment period $t$. The sum of generation over all sources must meet the demand in every demand block for every scenario (\ref{Demand}), and the output from each generation technology in every scenario for every demand block must be less than or equal to the available capacity in that investment period $t$ (\ref{CapMax}).

\begin{multline}
\label{Objective}
                \min_{x_{lt\sigma}^{g},y_{t\sigma}^{g}}{Z} \  =  \ \ \sum\limits_{t=1}^T \sum\limits_{g=1}^G \sum\limits_{\sigma=1}^\Sigma \sum\limits_{l=1}^L (OC_{t\sigma}^{g}x_{lt\sigma}^{g}H_{l})p(\sigma) + \\
                \sum\limits_{t=1}^T \sum\limits_{g=1}^G \sum\limits_{\sigma=1}^\Sigma (IC_{t\sigma}^{g}y_{t\sigma}^{g})p(\sigma) \ \ \   \; g \; \in \; G \;, \; t\;\in\;T\;, \; \sigma\;\in\;\Sigma, \; l \in L  
\end{multline}

 \begin{equation}\label{Demand}
               \sum\limits_{g=1}^G x_{lt\sigma}^{g}\ \  = D_{lt\sigma}, \; l \in L, \; t \in T, \; \sigma \in \Sigma \\
               \end{equation}
 \begin{equation}\label{CapMax}
               x_{lt\sigma}^g \leq y_{t\sigma}^g , \; t \in T, \; \sigma \in \Sigma, \; g \in G, \; l \in L
               \end{equation}

We reformulate the above SGEP as a Markov Decision Process with discrete time steps and continuous state and action spaces. The main components of the MDP correspond to the SGEP as follows:

\begin{itemize}
   \item Stages: The stages in the MDP directly correspond to the investment periods $t$ as formulated above. Each stage is a representative year, and the time interval between stage $t$ and stage $t+1$ may be multiple years.
   \item State Space: The state vector $s_t$ is defined as the current total capacity $z_{t\sigma}^{g}$ of each generation type $g$ in stage $t$ and scenario $\sigma$. 
    \item Action Space: The action vector in each state is the amount of new capacity $y_{t\sigma}^{g}$ to be built for each generation type $g$ for scenario $\sigma$ in stage $t$.  
    \item Exogenous Information: The exogenous uncertainty $\omega_{st}$ is the vector of uncertain operation and investment costs, which are elements of the cost parameters $OC_{t\sigma}^{g}$ and $IC_{t\sigma}^{g}$ in the objective function (\ref{Objective}).
    \item Transition Function: Determines the new total capacity available in stage $t$, after the addition of new capacity. Mathematically it will be as in (\ref{CapDiff}), with $z_{t\sigma}^g$ representing the generating capacity of type $g$ available in scenario $\sigma$, stage $t$.
     \begin{equation}\label{CapDiff}
               z_{t\sigma}^g=z_{{t-1}\sigma}^g+y_{t\sigma}^g , \; t \in T, \; \sigma \in \Sigma, \; g \in G
     \end{equation}

    \item Current Stage Cost: The current cost is the sum of the investment cost and the operating cost for stage $t$.
    
    \item Policy: The policy $\pi$ maps the the state variable $s_t \equiv z_{t\sigma}^g$ to the optimal investment decision $a_{st}^{*} = y_{t\sigma}^{g}$. The policy provides a decision rule for the optimal capacity investment plan for every state in every stage. 
\end{itemize}

The operations sub-problem in SGEP can include additional considerations and details such as including non-convex operational constraints of generators or representing transmission network flows and constraints, but at the expense of more computation time when obtaining each sample cost for state-action samples in each stage. Because our focus here is to compare the relative performance of alternative adaptive sampling approaches, we assume a simple form of the sub-problem, economic dispatch, which solves for the output of each generation type to meet demand in each of 16 load blocks that approximate the demand over an entire year. The stylized experiment for this work is based on publicly available data for the Western Electricity Coordinating Council (WECC) \cite{munoz2013engineering,ho2016planning,bukenberger2019approximate}. The load duration curve consists of sixteen blocks, four demand levels from the hourly observations for each of the four seasons. To account for renewable energy generation, which are not considered for investment, net demand is calculated for 16 load blocks. Current capacity for each generator type is also considered and is shown in Table \ref{tab: Current Capacity}. We neglect consideration of transmission constraints and generator operating constraints. The methods shown here are applicable to more detailed representations of the dispatch sub-problem, and other operations models could be substituted.

We explore the investment strategy for four generation technologies: Natural Gas Combustion Turbine (GT), Combined Cycle Natural Gas Turbine (CCGT), Coal and Nuclear. The model has three investment stages, each representing one year, with stages assumed to be 20 years apart. We do not consider construction time lags for new capacity; new capacity built in stage $t$ is available for generation. The uncertain parameters are natural gas prices and future carbon prices and are assumed to be uniformly distributed random variables between the upper and lower bounds provided in Table \ref{tab:Uncertainity bounds}.

\begin{table}[H]
  \begin{center}

    \begin{tabular}{| c | c | c | c | c |}
    \hline
    {} & \multicolumn{2}{c |}{\textbf{Gas Price}} & \multicolumn{2}{c |}{\textbf{Carbon Price}} \\
    \hline
  \textbf{Stage} & \textbf{Lower Bound} & \textbf{Upper Bound} & \textbf{Lower Bound} & \textbf{Upper Bound}\\ 
  \hline
  Stage 1 & 3.2 & 3.2 & 50 & 50 \\
  \hline
  Stage 2 & 3 & 7 & 0 & 100 \\
  \hline
  Stage 3 & 3 & 11 & 100 & 300\\
  \hline

    \end{tabular}
        \caption{Uncertainty bounds for natural gas price and carbon price}
    \label{tab:Uncertainity bounds}
  \end{center}
\end{table}

\begin{table}[H]
  \begin{center}

    \begin{tabular}{| c | c | c | c | c |}
    \hline
  \textbf{Fuel Type} & \textbf{Existing Capacity} \\ 
  \hline
  Natural Gas Combustion Turbine & 9760\\
  \hline
  Combined Cycle Gas Turbine & 12260\\
  \hline
  Coal & 9260 \\
  \hline
  Nuclear & 8260 \\
  \hline

    \end{tabular}
        \caption{Existing capacity for each generator type at $t=0$.}
    \label{tab: Current Capacity}
  \end{center}
\end{table}

\section{Numerical Results}
We compare the performance of the proposed QIS algorithm to two variants of epsilon greedy for the generation expansion problem described above. We assume a linear architecture for the approximate Q-function, but include bilinear and square basis functions to better capture the shape of the underlying value function. Choosing the best architecture for the approximate Q-function is an art in itself and is not the main focus of this paper; here we focus on methods to sample actions from a continuous space for a given architecture. Because the problem has a continuous feasible action space, many exploration algorithms for RL (e.g., R-Max, $E^3$, Knowledge gradient) cannot be applied without some discretization of the decision space.  We therefore leave comparisons to such methods to future work. The algorithms compared are:

\begin{enumerate}
    \item \textbf{QIS}: The proposed algorithm described in Section 5.
    
   \item \textbf{QIS-RE}: The  proposed variant of QIS that reevaluates the approximate Q-values, $\hat{q}$, for all prior samples and updates $qMax_{t}$ and $qMin_{t}$ less frequently.
    
    \item \textbf{Epsilon Greedy}: This algorithm prescribes the probability $\epsilon$ of randomly sampling the action; the optimal action for the current value function approximation, $\hat{q}$, for the sample state is selected with probability $1-\epsilon$ \cite{watkins1989learning,singh2000convergence}. We present the results for values of $\epsilon$ from 0 to 1 in increments of 0.1, held constant over all iterations.
    
    \item \textbf{Epsilon Decay}: This variant of epsilon greedy decreases the value of $\epsilon$ as the iterations progress. It therefore will have a higher proportion of exploration samples in early iterations and fewer exploration samples in later iterations. We assume that the value of $\epsilon$ is modified as in
    (\ref{Epsilon Decay Step}) with decay rate $\delta$. The parameter $\delta$ is derived from specified initial and terminal values $\epsilon_\textrm{Initial}$ and $\epsilon_\textrm{Final}$ respectively, and the the maximum number of iterations $K$ (\ref{Epsilon Decay}). 

\begin{equation}\label{Epsilon Decay Step}
\epsilon_\textrm{k} = \epsilon_\textrm{k-1}*\delta    
\end{equation} 

\begin{equation}\label{Epsilon Decay}
\epsilon_\textrm{Final} = \epsilon_\textrm{Initial}*\delta^\textrm{K}    
\end{equation} 

\item \textbf{SP}: A deterministic equivalent implementation of the multi-stage stochastic optimization model (\ref{Objective}-\ref{CapMax}) is used as a benchmark to validate the solutions of the various Q-learning experiments.  For the SP model, an exogenous scenario tree is constructed to represent the uncertainties in natural gas price and carbon price by normalizing the range of values to the interval [0,1] and discretizing in steps of 0.1 in each dimension. The investment decisions are continuous positive variables in the SP model. 
\end{enumerate}

The results for four algorithms (QIS, QIS-RE, Epsilon Greedy and Epsilon Decay) are presented below.  Default values for algorithm parameters assume $M=10$ samples per iteration, $K=900$ total iterations, with $\hat{K}$=20 for QIS-RE. Q-function approximation, $\hat{q_{t}}$, coefficients are initialized to zero, and $qMin_{t}$ and $qMax_{t}$ are equal to 0 and 1, respectively, unless otherwise specified. All experiments assume a Temporal Difference learning rate $\lambda = 0.1$ to update the approximation coefficients $\theta$ (\ref{Update Q Value}).  The models using QIS, QIS-RE, Epsilon Greedy and Epsilon Decay were implemented in Matlab. The non-linear program solved to obtain exploitation samples in Epsilon Greedy and Epsilon Decay were solved using the interior point algorithm with the fmincon function in Matlab. The simulation of the policy in the Q-learning algorithms to obtain the total cost was performed on the same scenario tree used in the deterministic equivalent implementation. The SP model and simulation of the Q learning policy are implemented in GAMS and solved by CPLEX 12.8.0. All experiments are performed on a high-performance computing cluster with 2.2 GHz Intel Xeon Processor, 24 cores, 128 GB RAM, and 40 Gbps Ethernet.

\subsection{Stability of the optimal solution}\mbox{}
An effective sampling algorithm within ADP should facilitate convergence to a policy that is close to optimal in a finite number of iterations.  Moreover, the solution should be robust to the algorithmic parameters, and should be consistent across multiple independent repetitions.  We evaluate the consistency of the policy obtained from 10 replications of each algorithm. The QIS and QIS-RE algorithms obtain the same first stage policies for all replications, and exhibit variation in the optimal cost of less than 0.1\% (Table \ref{tab:Result_Compare}). We also show the solution from the SP model for comparison.

\begin{figure}[H]
\centering
   \includegraphics[width=1\linewidth]{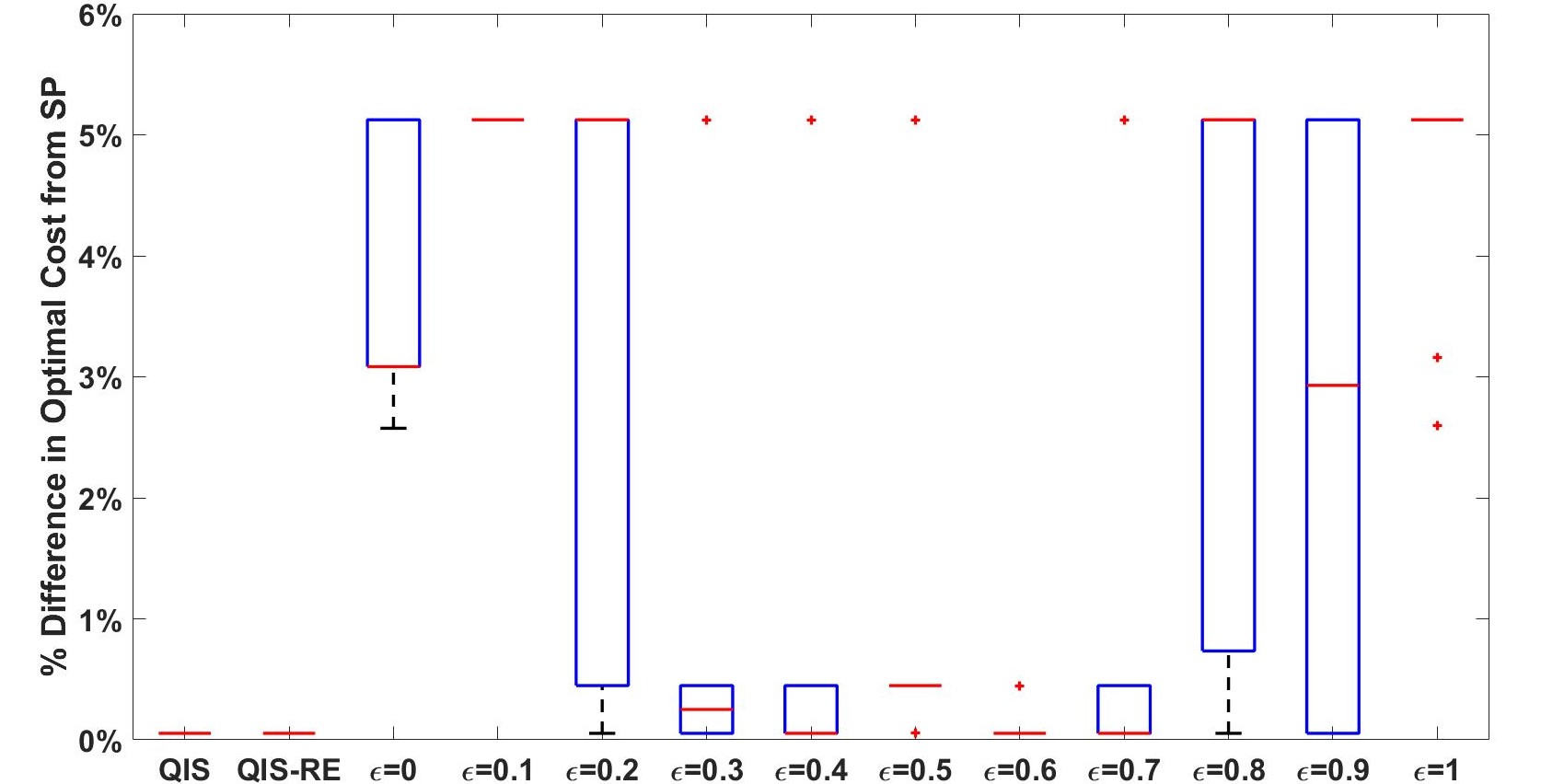}
   \caption{Comparison of \% difference in optimal expected cost relative to the SP solution from simulating the optimal policies obtained in Epsilon Greedy for $\epsilon$ values from 0 (pure exploit) to 1 (pure explore) in increments of 0.1 for 10 replications with QIS and QIS-RE. Boxes enclose the 50\% interval, midlines indicate the median value, whiskers indicate the 90\% interval, and outlies are shown as '+'.}
  \label{fig:Different Epsilon}
\end{figure}

In Epsilon Greedy, the critical algorithmic parameter is $\epsilon$, the probability of randomly sampling the next action, as opposed to choosing the optimal action with respect to the current approximate Q-function $\hat{q}$. In general, this value is chosen heuristically and is application-specific. Figure \ref{fig:Different Epsilon} shows the \% difference in the expected optimal cost from the Epsilon Greedy version relative to the benchmark optimal cost from the Stochastic Programming version across 10 replications by simulating the optimal policies for epsilon values ranging from 0 to 1 in increments of 0.1. The choice of 0.6 for $\epsilon$ appears to provide the best results, with most repetitions diverging from the SP solution by less than 0.1\% except for a few outliers. For values of $\epsilon$ between 0.0 and 0.3, the infrequent exploration of actions often leads to premature convergence to a sub-optimal solution. For values of $\epsilon$ above 0.7, the frequent exploration results in greater variance across different repetitions.

\begin{figure}[H]
   \includegraphics[width=1\linewidth]{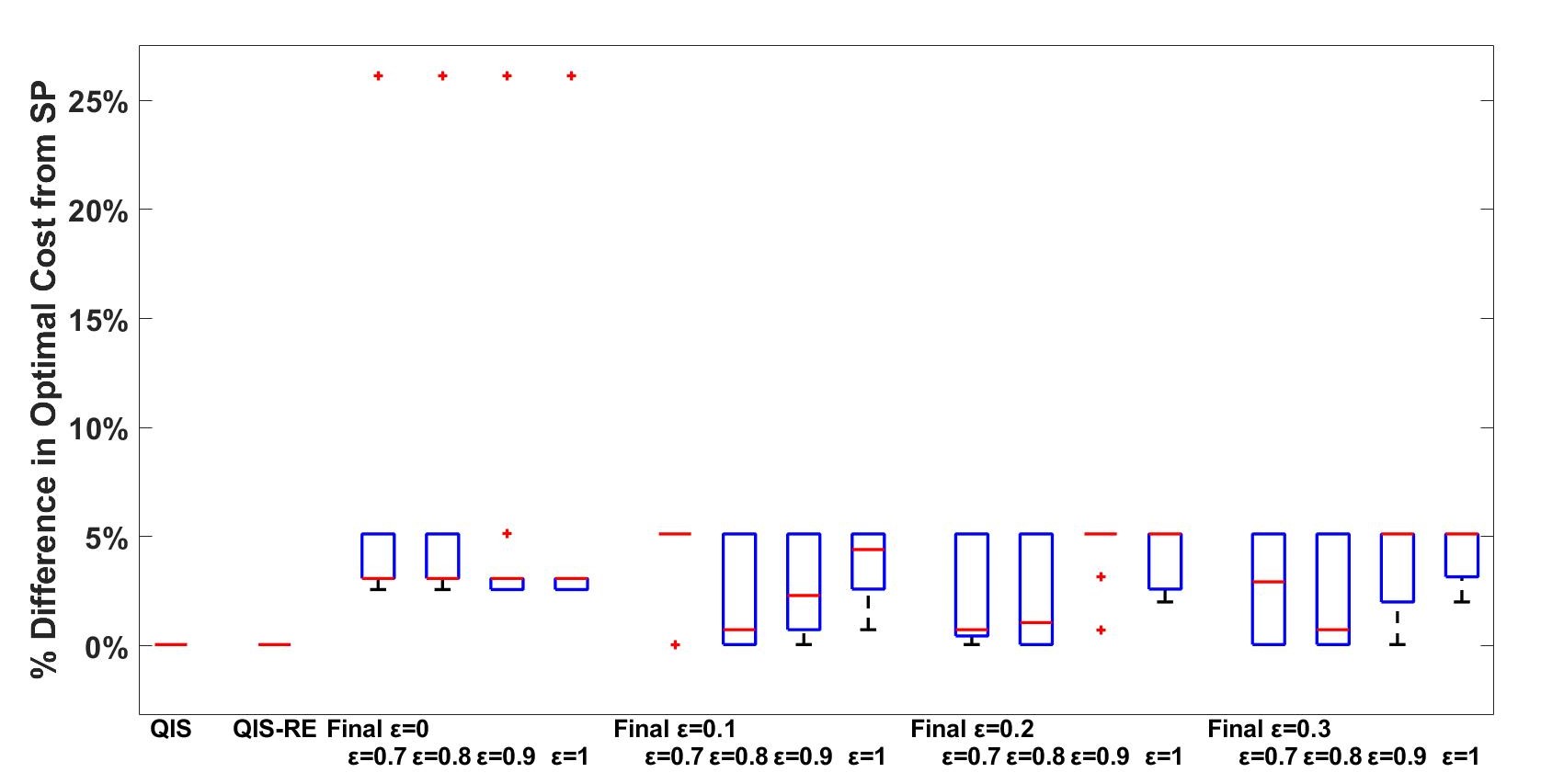}
  \caption{Comparison of \% difference in optimal expected cost relative to the SP solution from simulating the optimal policies obtained in Epsilon Decay for different starting and final $\epsilon$ values with reductions according to Equation \ref{Epsilon Decay} for 10 replications with QIS and QIS-RE. Boxes enclose the 50\% interval, midlines indicate the median value, whiskers indicate the 90\% interval, and outlies are shown as '+'.}
  \label{fig:Decay_Diff_Epsilon}
\end{figure}

Epsilon Decay is a variant of Epsilon Greedy in which the exploration rate decreases as a function of the iteration count.  The value of $\epsilon$ decreases at an exogenously specified rate; as more cumulative sample costs are observed, actions are randomly sampled less frequently. We parameterize the decay rate as in (\ref{Epsilon Decay}), where the user specifies the initial and final values of $\epsilon$. We present the results from the Epsilon Decay version for all combinations of final values of $\epsilon$ \{0,0.1,0.2,0.3\} and initial values \{0.7,0.8,0.9,1\} (Figure \ref{fig:Decay_Diff_Epsilon}). For most choices of the initial, final, and decay rate, the Epsilon Decay version can sometimes converge to a solution that is close to the benchmark optimal cost, but the variance across independent repetitions of the algorithm is large in most cases.

\begin{table}[H]
  \begin{center}

    \begin{tabular}{| c | c | c | c | c | c | c |}
    \hline
  \textbf{Algorithm} & \textbf{\% Nuclear} & \textbf{\% Coal} & \textbf{\% CCGT} & \textbf{\% GT} & \textbf{Total Cost} & \textbf{\% Difference from SP} \\ 
  \hline
  QIS & 33.3 & 0.0 & 33.3 & 33.3 & 2.565E+11 & 0.05\% \\
  \hline
  QIS-RE & 33.3 & 0.0 & 33.3 & 33.3 & 2.565E+11 & 0.05\% \\
  \hline
  Epsilon Greedy & 50.0 & 0.0 & 0.0 & 50.0 & 2.575E+11& 0.45\%\\
  \hline
  Epsilon Decay  &  0.0 & 0.0 & 0.0 & 100.0 & 2.695E+11 & 5.12\% \\
  \hline
  SP & 32.0 & 0.0 & 32.0 & 36.0 & 2.564E+11 & \\
  \hline

    \end{tabular}
        \caption{Comparison of the first stage policy obtained for QIS, QIS-RE,  Epsilon Greedy with $\epsilon=0.5$, Epsilon Decay with $\epsilon_{Initial}=0.7$ and $\epsilon_{Final}=0.2$ and SP for one of the replications of the algorithm. \% Nuclear, \% Coal, \% CCGT and \% GT represent the percentage of the first stage capacity decisions. Total Cost is obtained by stochastic simulation of the policy for the same scenario tree as SP. \% Difference from SP refers to the \% difference in optimal expected cost relative to the SP solution from simulating the optimal policies obtained in QIS, QIS-RE, Epsilon Greedy, Epsilon Decay for one of the replications of the algorithm.}
    \label{tab:Result_Compare}
  \end{center}
\end{table}


\subsection{Balancing Exploration and Exploitation }\mbox{}
Sampling strategies for ADP must balance exploration, sampling actions to improve the approximate value function, with exploitation, sampling the optimal action to improve the estimate of its expected value. Epsilon Greedy achieves this balance with the parameter $\epsilon$, the probability of sampling the next action randomly. QIS balances exploration and exploitation using the qRatio, which accepts proposed sample actions proportionally to their current estimated approximate Q-Value, $\hat{q}$. The QIS algorithm has a non-zero probability of sampling any feasible action, but in contrast to Epsilon Greedy, does not assume a uniform distribution over the feasible action space. By sampling proportional to the approximate Q-function, QIS will continue to explore more actions in the neighborhood of the current approximation optimum. Samples from Epsilon Greedy form a mixture distribution of uniformly random samples (explore) and a single value at the current optimum (exploit). 

\begin{figure}[H]
\centering
\begin{subfigure}{.475\textwidth}
  \centering
  \includegraphics[width=\textwidth]{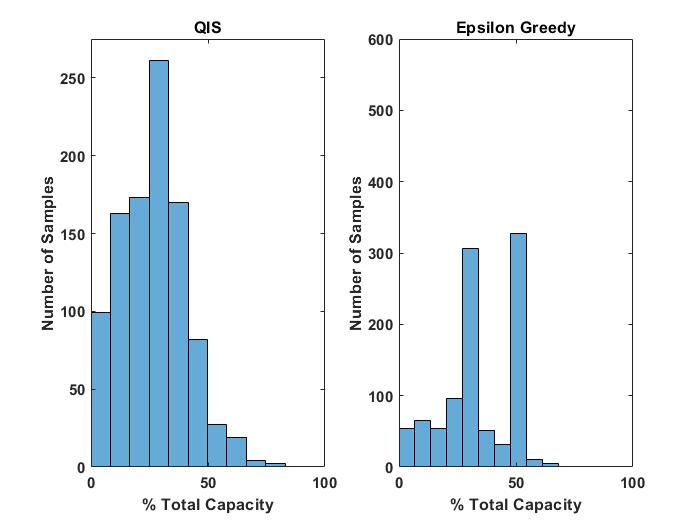}\hfill
  \caption{GT Samples}
  \label{fig:IS Greedy GT}
\end{subfigure}%
\begin{subfigure}{.475\textwidth}
  \centering
  \includegraphics[width=\textwidth]{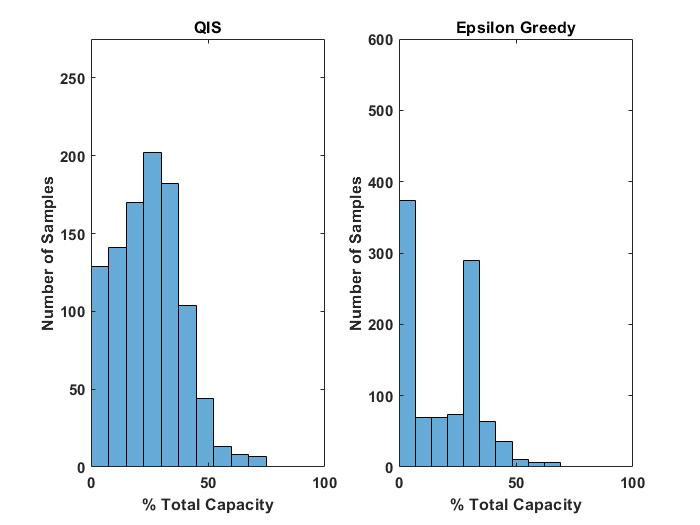}\hfill
  \caption{CCGT Samples}
  \label{fig:IS Greedy CCGT }
\end{subfigure}
\vskip\baselineskip
\begin{subfigure}{.475\textwidth}
  \centering
  \includegraphics[width=\textwidth]{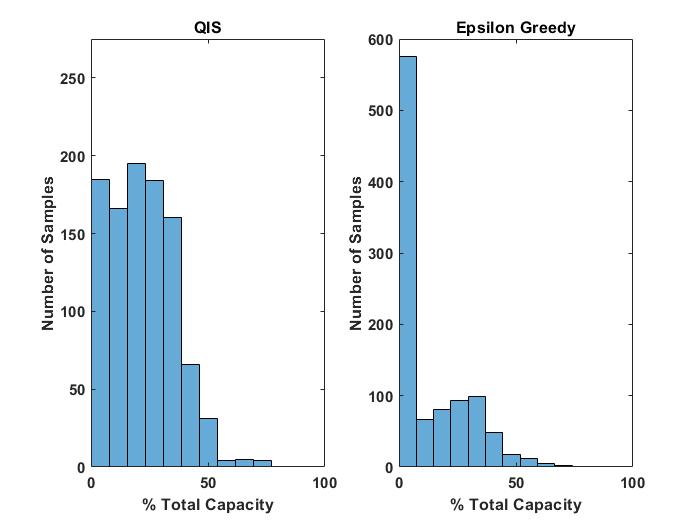}\hfill
  \caption{Coal Samples}
  \label{fig:IS Greedy Coal }
\end{subfigure}
\begin{subfigure}{.475\textwidth}
  \centering
  \includegraphics[width=\textwidth]{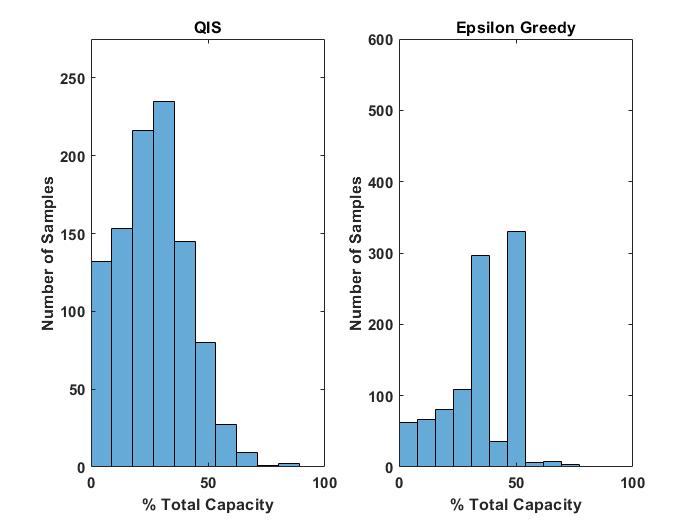}\hfill
  \caption{Nuclear Samples}
  \label{fig:IS Greedy Nuclear}
\end{subfigure}
\caption{Final 1000 action samples for QIS and Epsilon Greedy as a percentage of Stage 1 new capacity for Natural Gas Combustion Turbine (\ref{fig:IS Greedy GT}), Combined Cycle Gas Turbine (\ref{fig:IS Greedy CCGT }), Coal (\ref{fig:IS Greedy Coal }) and Nuclear (\ref{fig:IS Greedy Nuclear}) technologies for one representative replication. Epsilon Decay and QIS-RE have not been shown because of similar behavior to Epsilon Greedy and QIS respectively. Note the change in the $Y$-axis scale between QIS and Epsilon Greedy.}
\label{fig: IS Greedy Samples}
\end{figure}

The difference between sampling distributions from QIS and Epsilon Greedy is illustrated in Figure \ref{fig: IS Greedy Samples}, which shows the histogram of the final 1000 samples for new capacity of each generator technology in Stage 1 for one representative replication. Sample actions are expressed as a percentage of the required new capacity for each generator technology; the percentage of total new capacity built must sum to $100\%$. The samples from QIS for new Nuclear, CCGT, and GT capacity additions more frequently fall close to the optimal decisions from the SP (Table \ref{tab:Result_Compare}). Samples of new coal capacity additions shift towards nearly equally spaced samples for capacity below 40\%. Epsilon Greedy can sometimes result in frequently sampled actions from two distinct regions of the decision space (Figure \ref{fig: IS Greedy Samples}).  In this example, sampled actions oscillate between the optimal SP new capacity shares (Table \ref{tab:Result_Compare}) and a suboptimal policy of 50\% Nuclear and 50\% GT. 

The use of temporal differencing (\ref{Update Q Value}) to update the approximation for continuous action spaces limits the incremental information value of repeatedly visiting the exact same action. Early iterations can often lead to an approximation with a local optimum.  Sampling very suboptimal actions may not be sufficient to move the approximation away from the local optimum towards a global optimum.  A sampling strategy that more frequently explores the immediate neighborhood of the optimal action, as in QIS, is better able to avoid getting stuck on a suboptimal solution.  All replications performed have behavior similar to the one presented in Figure \ref{fig: IS Greedy Samples}.

\subsection{Sensitivity to the number of iterations and the number of samples per iteration}
\begin{figure}[H]
\centering
\begin{subfigure}{.3\textwidth}
  \centering
  \includegraphics[width=\textwidth]{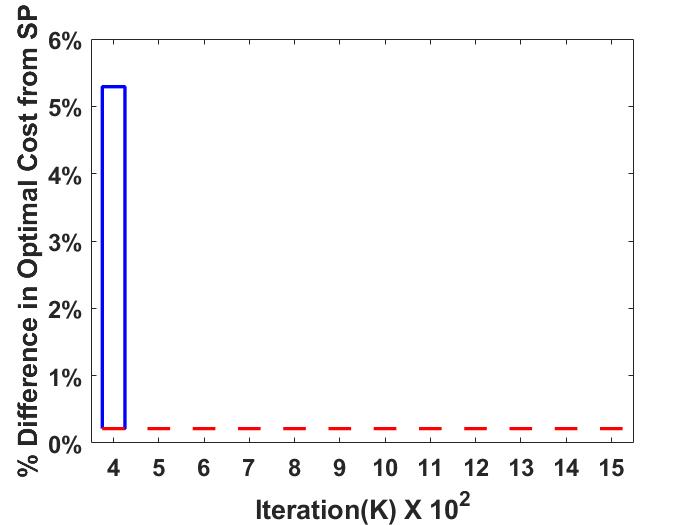}\hfill
  \caption{QIS \& QIS-RE}
  \label{fig:Iter_IS}
\end{subfigure}%
\begin{subfigure}{.3\textwidth}
  \centering
  \includegraphics[width=\textwidth]{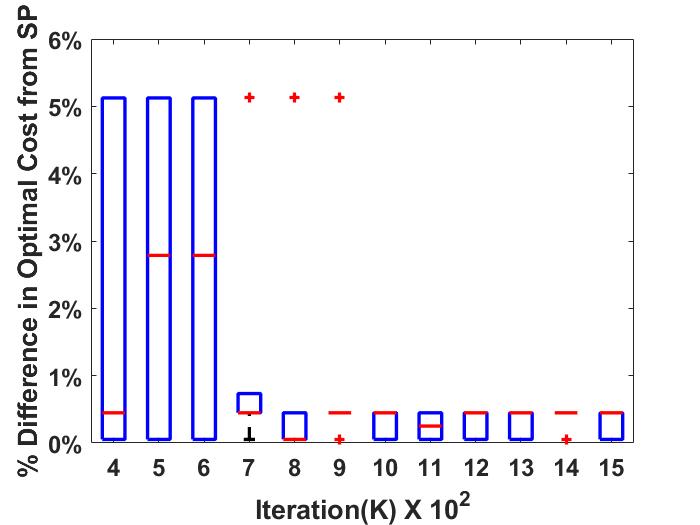}\hfill
  \caption{Epsilon Greedy}
  \label{fig: Iter_greedy}
\end{subfigure}
\begin{subfigure}{.3\textwidth}
  \centering
  \includegraphics[width=\textwidth]{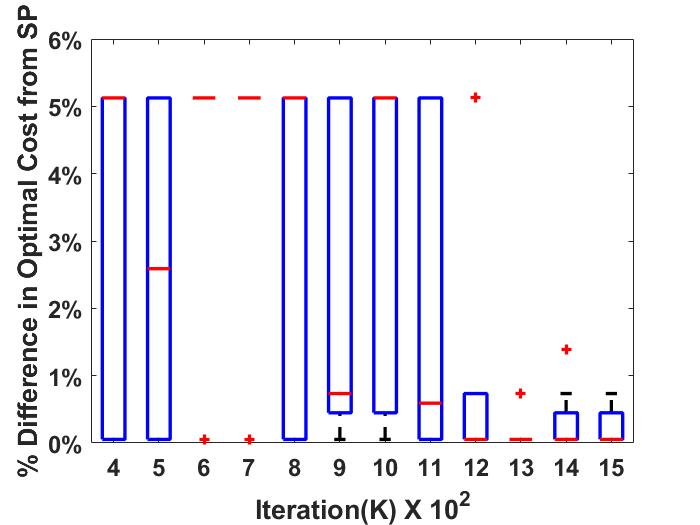}
  \caption{Epsilon Decay with $\epsilon_{Initial}$=0.7 and $\epsilon_{Final}$=0.2}
  \label{fig:Iter_Decay}
\end{subfigure}
\caption{Comparison of \% difference in optimal expected cost relative to the SP solution from simulating the optimal policies obtained by QIS and QIS-RE (\ref{fig:Iter_IS}), Epsilon Greedy with $\epsilon=0.5$ (\ref{fig: Iter_greedy}) and Epsilon Decay (\ref{fig:Iter_Decay}) with $\epsilon_{Initial}=0.7$ and $\epsilon_{Final}=0.2$. Shown are results for $K=400$ to $K=1500$ in increments of 100 for 10 replications. Boxes enclose the 50\% interval, midlines indicate the median value, whiskers indicate the 90\% interval, and outlies are shown as '+'.}
\label{Iteration Count}
\end{figure}

\mbox{}
In this section, we compare the rate of convergence across multiple replications to the optimal policy for QIS, QIS-RE, Epsilon Greedy and Epsilon Decay. Figure 5 represents the \% difference in the estimated expected optimal cost of the policy, relative to the expected cost from the SP model, for 10 replications after different numbers of iterations. We show results for one instance of Epsilon Greedy that assumes $\epsilon = 0.5$ (Fig. \ref{fig: Iter_greedy}), and one instance of Epsilon Decay that assumes $\epsilon_{Initial} = 0.7$ and $\epsilon_{Final} = 0.2$ (Fig. \ref{fig:Iter_Decay}). QIS and QIS-RE both consistently converge within 500 iterations to the optimal SP policy for all replications.  The other two methods converge after 1000 iterations, but exhibit higher variability across replications. For this example, Epsilon Greedy finds better solutions with less variability than Epsilon Decay. One possible reason for the poor solution quality from Epsilon Decay in these results is that it shifts to more exploitation samples before the approximation is sufficiently accurate. 

\begin{figure}[H]
\centering
\includegraphics[width=.8\textwidth]{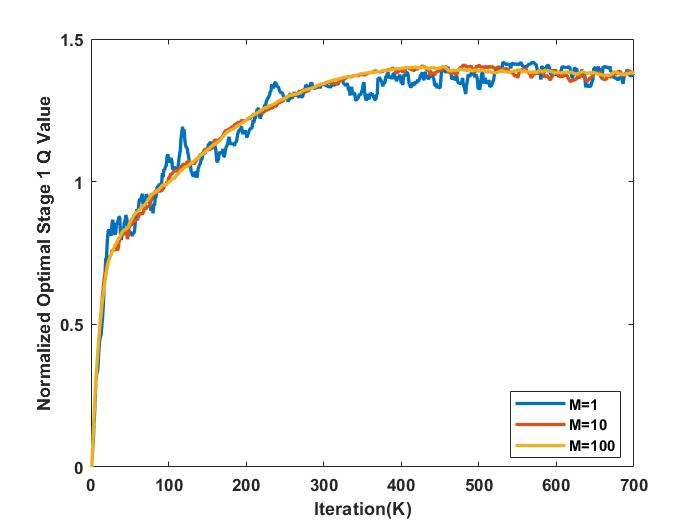}\hfill
\caption{Normalized Optimal Stage 1 Q-value from the QIS algorithm for iterations $K=1$ to $700$ in increments of $1$ for samples size $M=1, 10$, and $100$ for one representative replication.}
\label{fig: Line_Plot}
\end{figure}

We also explore the effect of different sample sizes within each iteration. We repeat the experiments with $M = 1, 2, 5, 10, 25, 50 \text{ and } 100$, and present the \% difference in the estimated expected optimal cost of the policy, relative to the optimal expected cost from the SP solution.  The results are provided at two iteration counts, $K = 500 \text{ and } 900$ in Table \ref{tab: Samples Sensitivity}. For small sample sizes of 1, 2 and 5, the QIS results exhibit high variability at $K=500$, but all sample sizes converge to the same policy by $K=900$ iterations, although some outliers are still observed for sample sizes of $M = 1$ and $2$ for $K \geq 900$ (not shown). The QIS-RE algorithm exhibits behavior similar to that of QIS, but with greater variability in expected optimal costs for small samples sizes.  Small sample sizes and infrequent re-evaluation of Q-values can results in poor estimates of the gradient and the temporal difference error in Equation (\ref{Update Q Value}). Larger per-iteration sample sizes enable better estimates of the gradient and error, which improves the solution quality and stability. The impact of the per-iteration sample size for the QIS algorithm is illustrated in Figure \ref{fig: Line_Plot} in terms of the normalized Stage 1 optimal Q-function value; larger sample sizes result in a smoother evolution. 

The estimate expected optimal costs for Epsilon Greedy with $\epsilon = 0.5$ have significant variability for all sample sizes at $K=500$ iterations. After $K = 900$ iterations, the variance in the optimal cost across replications is smaller for sample sizes of 2 or more. However, the variability of solutions when using Epsilon Greedy is still greater than the solutions when using QIS or QIS-RE. The solutions using Epsilon Decay exhibit the greatest variation across replications for all sample sizes, relative to the other sampling methods. Further, Epsilon Decay sometimes finds an optimal cost close to the true optimum after $K = 500$ iterations, but then diverges from this optimum with additional iterations. This suggests that as the frequency of exploitation increases, the model is choosing sub-optimal decisions because of local optima when solving for the optimal action or because of insufficient samples of other feasible actions.

\begin{table}[H]
\resizebox{\columnwidth}{!}{%
\centering

    \begin{tabular}{| c | c | c | c | c | c | c | c | c | c | c | c | c | c | c|}
    \hline
    \multicolumn{1}{|l}{} &  & \multicolumn{3}{c |}{\textbf{QIS}} & \multicolumn{3}{c |}{\textbf{QIS-RE}} & \multicolumn{3}{c |}{\textbf{Epsilon Greedy}}& \multicolumn{3}{c |}{\textbf{Epsilon Decay}} \\
    \hline
  \textbf{K} &\textbf{M} & \textbf{Min(\%)} & \textbf{Med(\%)} & \textbf{Max(\%)} & \textbf{Min(\%)} & \textbf{Med(\%)} & \textbf{Max(\%)} & \textbf{Min(\%)} & \textbf{Med(\%)} & \textbf{Max(\%)} & \textbf{Min(\%)} & \textbf{Med(\%)} & \textbf{Max(\%)}\\ 
  \hline
   \multirow{7}{*}{500} & \textbf{1} & 0.05 & 0.05 & 5.12 & 0.05 & 0.05 & 5.12 & 0.05 & 5.12 & 5.12 & 0.05 & 0.8 & 5.12 \\\cline{2-14}
     & \textbf{2} &     0.05 & 0.05 & 5.12 & 0.05 & 0.05 & 5.12 & 0.05 & 5.12 & 5.12 & 0.05 & 5.12 & 5.12 \\\cline{2-14}
      & \textbf{5} &    0.05 & 0.05 & 5.12 & 0.05 & 0.05 & 5.12 & 0.05 & 5.12 & 5.12 & 0.05 & 5.12 & 5.12 \\\cline{2-14}
      & \textbf{10} &    0.05 & 0.05 & 0.05 & 0.05 & 0.05 & 0.05 & 0.05 & 2.79 & 5.12 & 0.05 & 2.59 & 5.12 \\\cline{2-14}
     & \textbf{25} &     0.05 & 0.05 & 0.05 & 0.05 & 0.05 & 0.05 & 0.05 & 2.59 & 5.12 & 0.05 & 0.05 & 5.12 \\\cline{2-14}
    & \textbf{50} &      0.05 & 0.05 & 0.05 & 0.05 & 0.05 & 0.05 & 0.05 & 5.12 & 5.12 & 0.05 & 0.05 & 5.12 \\\cline{2-14}
     & \textbf{100} &     0.05 & 0.05 & 0.05 & 0.05 & 0.05 & 0.05 & 0.05 & 5.12 & 5.12 & 0.05 & 0.05 & 5.12 \\ \hline

\multirow{7}{*}{900} &\textbf{1} &    0.05 & 0.05 & 0.05 & 0.05 & 0.05 & 5.12 & 0.45 & 0.69 & 5.12 & 0.09 & 0.45 & 5.12 \\\cline{2-14}
    {} & \textbf{2} &   0.05 & 0.05 & 0.05 & 0.05 & 0.05 & 0.05 & 0.05 & 0.45 & 0.45 & 0.05 & 0.45 & 5.12 \\\cline{2-14}
    {} & \textbf{5} &    0.05 & 0.05 & 0.05 & 0.05 & 0.05 & 0.05 & 0.05 & 0.25 & 0.45 & 0.73 & 5.12 & 5.12 \\\cline{2-14}
    {} & \textbf{10} &    0.05 & 0.05 & 0.05 & 0.05 & 0.05 & 0.05 & 0.05 & 0.45 & 5.12 & 0.05 & 0.73 & 5.12 \\\cline{2-14}
    {} & \textbf{25} &    0.05 & 0.05 & 0.05 & 0.05 & 0.05 & 0.05 & 0.05 & 0.05 & 5.12 & 0.05 & 5.12 & 5.12 \\\cline{2-14}
    {} & \textbf{50} &    0.05 & 0.05 & 0.05 & 0.05 & 0.05 & 0.05 & 0.05 & 0.05 & 0.45 & 0.05 & 2.93 & 5.12 \\\cline{2-14}
    {} & \textbf{100} &    0.05 & 0.05 & 0.05 & 0.05 & 0.05 & 0.05 & 0.05 & 0.05 & 0.45 & 0.05 & 2.93 & 5.12 \\ \hline
   
    \end{tabular}%
    }
        \caption{Minimum, Median and Maximum \% difference in the optimal expected cost relative to the SP solution from simulating the optimal policies obtained by QIS, QIS-RE, Epsilon Greedy with $\epsilon=0.5$, and Epsilon Decay with $\epsilon_{Initial}=0.7$ and $\epsilon_{Final}=0.2$ for $K=500$ and $900$ with sample size $M=1, 2, 5, 10, 25, 50$, and $100$ for $10$ replications.}
    \label{tab: Samples Sensitivity}

\end{table}

\subsection{Computational Effort}
\begin{figure}[H]
\centering
\includegraphics[width=0.8\textwidth]{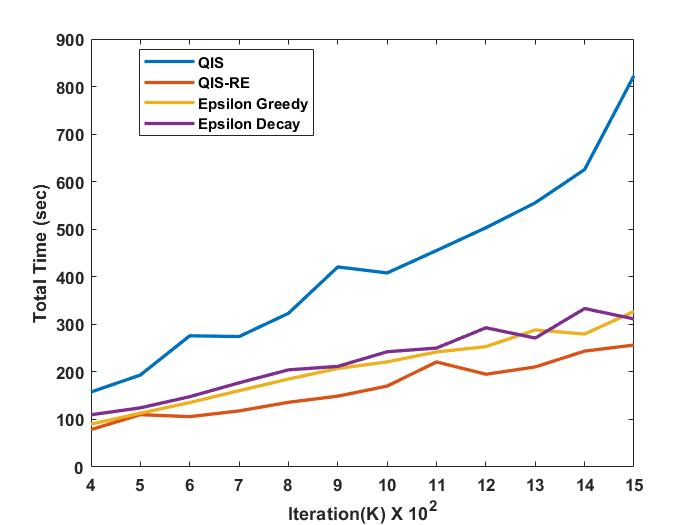}\hfill
\caption{Total computation time (sec), averaged across 10 replications, for QIS, QIS-RE with $\hat{K}=20$, Epsilon Greedy with $\epsilon=0.5$ and Epsilon Decay with $\epsilon_{Initial}=0.7$ and $\epsilon_{Final}=0.2$ for iterations $K=400$ to $1500$ in increments of $100$ with $M=10$.}
\label{fig: Time_Plot}
\end{figure}

In addition to producing high quality, low variance estimates of the optimal cost and policy, an effective sampling method should also have minimal computation time.  Less computation time expended obtaining samples allows for larger problems to be solved.  We compare the total computation time for all four sampling methods, averaged across 10 replications, as a function of iteration count (Figure \ref{fig: Time_Plot}). The QIS algorithm requires the most computation time, followed by Epsilon Decay, Epsilon Greedy and QIS-RE, which requires the least computation time. To understand the differences across methods, we distinguish the time required within three distinct algorithmic sub-tasks (Table \ref{tab:Computational Time}): the "Sampling Time" is the time required to obtain the sample action for the given state; the "Evaluation Time" is the time required to recompute all approximate Q-values and update the estimated $qMax_{t}$ and $qMin_{t}$; and "Other Time" includes all other tasks including updating the approximation coefficients and parallel computing communication. The Evaluation Time is not required for Epsilon Greedy and Epsilon Decay, because these methods do not update estimates of $qMax_{t}$ and $qMin_{t}$. QIS requires the most computational effort because of the need to recompute the approximate Q-values for all samples from previous iterations after every iteration. Epsilon Greedy and Epsilon Decay require similar computation time; Epsilon Decay requires slightly more time because the increasing frequency of exploit samples requires more optimizations of the Q-function to obtain these sample actions. Both QIS and QIS-RE only evaluate the Q-function, but do not need to optimize it. QIS-RE requires the least total computation time among these algorithms because it does not have to optimize the value function, as in Epsilon Greedy and Epsilon Decay, and because it reduces the frequency of recomputing all Q-values, which takes the majority of the time for the QIS algorithm. In the example shown in Table \ref{tab:Computational Time}, QIS-RE recomputes the approximate Q-values every 20 iterations ($\hat{K}=20$), which significantly reduces the total computation time.

The QIS-RE does require an additional algorithmic parameter, $\hat{K}$, the frequency of recomputing $qMax_{t}$ and $qMin_{t}$. Table \ref{tab: Evaluation Sensitivity} presents the proportion of the replications that converge to the optimal Stage 1 policy (as in Table \ref{tab:Result_Compare}) for several values of $\hat{K}$, at different iteration counts, and for two different sample sizes. If QIS-RE never updates the Q-values, convergence to the optimal policy requires more iterations. Relatively infrequent reevaluation, such as every 50 or 100 iterations, has nearly the same performance as more frequent reevaluation, but requires less computation time. Furthermore, larger sample sizes ($M=10$) provide better solution consistency independent of the frequency of reevaluation.  Even the No Update version converges to the optimal policy with less variance than either Epsilon Greedy or Epsilon Decay if larger per-iteration sample sizes are used. We recommend updating the estimates of $qMax_{t}$ and $qMin_{t}$, even if infrequently, because it accelerates convergence.

\begin{table}[H]
  \begin{center}

    \begin{tabular}{| c | c | c | c | c | c | c |}
    \hline
  \textbf{Algorithm} & \textbf{Sampling Time (sec)} & \textbf{Evaluation Time (sec)} & \textbf{Other Time (sec)} &\textbf{Total Time}\\ 
  \hline
  QIS & 58.7 & 267.7 & 94.4 & 420.8\\
  \hline
  QIS-RE & 53.1 & 12.8 & 83.1 & 149.0 \\
  \hline
  Epsilon Greedy & 102.9 & - & 104.3 & 207.2 \\
  \hline
  Epsilon Decay & 113.5 & - & 98.1 & 211.6 \\
  \hline
 
    \end{tabular}
        \caption{Computation time (sec) for Sampling Time, Evaluation Time, Other Time, and Total Time, averaged across 10 replications of QIS, QIS-RE with $\hat{K}=20$, Epsilon Greedy with $\epsilon=0.5$ and Epsilon Decay with $\epsilon_{Initial}=0.7$ and $\epsilon_{Final}=0.2$. Results are for $K=900$ iterations and $M=10$ samples.}
    \label{tab:Computational Time}
  \end{center}
\end{table}

\begin{table}[H]
\centering
\begin{center}
    \begin{tabular}{| c | c | c | c | c | c | c | c | c | c | c |}
    \hline
    {} & \multicolumn{4}{c |}{\textbf{M=1}} & \multicolumn{4}{c |}{\textbf{M=10}} \\
    \hline
  \textbf{$\hat{K}$ vs K} & \textbf{K=400} & \textbf{K=500} & \textbf{K=600}& \textbf{K=700} & \textbf{K=400} & \textbf{K=500} & \textbf{K=600}& \textbf{K=700}\\ 
  \hline
  $\hat{K}$=1 & 40\% & 70\% & 90\% & 90\% & 60\% & 100\% & 100\% & 100\%  \\
  \hline
  $\hat{K}$=10 & 50\% & 100\% & 100\% & 90\% & 50\% & 100\% & 100\% & 100\%   \\ 
  \hline
  $\hat{K}$=20 & 60\% & 80\% & 90\% & 100\% & 60\% & 100\% & 100\% & 100\%   \\
  \hline
  $\hat{K}$=50 & 80\% & 100\% & 90\% & 90\% & 50\% & 100\% & 100\% & 100\%   \\
  \hline
  $\hat{K}$=100 & 90\% & 60\% & 90\% & 100\% & 70\% & 100\% & 100\% & 100\%  \\
  \hline
  No Update & 60\% & 60\% & 60\% & 70\% & 30\% & 90\% & 90\% & 100\%  \\
  \hline
    \end{tabular}
        \caption{Frequency within 10 replications that the optimal Stage 1 policy from QIS-RE is the policy shown in Table \ref{tab:Result_Compare}. Results show current optimal policy for different reevaluation frequencies $\hat{K}$, for sample sizes of $M=1$ and $10$, and for iteration counts $K=400, 500, 600$, and $700$.}
    \label{tab: Evaluation Sensitivity}
  \end{center}
\end{table}

\section{Conclusions}
We have presented the QIS algorithm, a novel scheme for sampling continuous action spaces in Approximate Dynamic Programming formulations. The approach uses importance sampling to sample actions proportional to the current state-action value function approximation. In contrast to methods that rely on selecting the current optimal action or sampling uniformly across the feasible action space, we have demonstrated superior solution quality and lower variance across replications when using QIS.  One advantage of this algorithm is the absence of tuning parameters when sampling actions, because it relies only on the approximate value function.  Another advantage is that it avoids having to solve a (possibly nonlinear) optimization to select the optimal action within the sampling step. One disadvantage of the QIS algorithm is the additional computational effort required to reevaluate prior samples after every iteration.  We have also proposed a variant on the algorithm, QIS-RE, that requires less computational effort and achieves the same solution quality and stability across replications as QIS.   QIS-RE does require an additional parameter for the frequency of reevaluations, although our results show that convergence is relatively insensitive to different assumptions.

Both the QIS and the QIS-RE algorithms can scale to larger multi-dimensional continuous state and action spaces, and they can be applied to any value function approximation architecture. The action space in the case study presented was continuous but this method is likely applicable to discrete decision spaces as well; this is left for future work. A potential improvement to the algorithm, not explored here, would be to gradually increase the relative density of sampled actions very close to the optimal action as iterations progress by means of an exponential parameter applied to the qRatio. This parameter would have low values in the initial iterations to promote exploration (more samples accepted farther from the optimal action) and higher values later to increase exploitation (most accepted samples very close to optimal action). The approximation architecture used in this paper for the Q-function was linear, but we expect this sampling scheme to be applicable to non-linear function approximation architectures, such as neural networks. The QIS sampling approach could be especially advantageous for continuous decision spaces when Epsilon Greedy methods are difficult to apply within a neural network based approximation. These areas of exploration are left for future work.

Much of the work on the application presented here, generation capacity expansion for electric power systems, has relied on methods from within the stochastic programming paradigm.   As multi-stage stochastic problems grow large, the decomposition methods within that paradigm often rely on dual variables and other assumptions about the functional form (e.g., convexity) of the system being studied.  Many of the applied questions in the area of energy transition and achievement of a low carbon energy system will require larger multi-stage problems that cannot omit critical non-convex features of the power system, which can limit the applicability of traditional decomposition schemes.   Although the case study in this work was intentionally simplified to provide an optimal benchmark, the method here offers new avenues for solving large multi-stage stochastic problems with non-convex features. An exploration of these problems is another area ripe for future work.

\section{Acknowledgments}
The authors gratefully acknowledge Dr. Uday Shanbhag, Dr. Eunhye Song, and Kshitij Dawar for helpful feedback and comments on this work. This research is supported by the U.S. Department of Energy, Office of Science, Biological and Environmental Research Program, Earth and Environmental Systems Modeling, MultiSector Dynamics, Contract No. DE-SC0016162. The opinions, findings, and conclusions or recommendations expressed in this material are those of the author(s) and do not necessarily reflect the views of the US Department of Energy. Computational experiments for this research were performed on the Pennsylvania State University's Computational and Data Sciences Roar supercomputer.  

\printbibliography

\end{document}